\documentclass[10pt,a4paper,leqno,english]{article}
\usepackage{graphicx}
\usepackage[T1]{fontenc}
\usepackage{babel}
\newtheorem{example}{Example}[section]

\newtheorem{proposition}[example]{Proposition}
\newtheorem{theorem}[example]{Theorem}

\newtheorem{rem}[example]{Remark}

\usepackage{enumerate,amsmath,amssymb,latexsym}

\newenvironment{proof}[1][Proof]{\noindent\textit{#1.} }{\ \rule{0.5em}{0.5em}}
\date{}
\newenvironment{keywords}[1][Keywords]{\noindent\textbf{#1.} }{\ \rule{0.5em}{0.5em}}
\newenvironment{M}[1][M.S.C. 2000]{\noindent\textbf{#1.} }{\ \rule{0.5em}{0.5em}}
\newenvironment{abs}[1][Abstract]{\noindent\textbf{#1.} }{\ \rule{0.5em}{0.5em}}

\begin{document}

\title{Biharmonic curves on \textit{LP}-Sasakian
manifolds }
\author{Sad\i k Kele\c{s}, Selcen Y\"{u}ksel Perkta\c{s}, Erol K\i l\i
\c{c}} \maketitle

\abs{In this paper we give necessary and sufficient conditions for
spacelike and timelike curves in a conformally flat, quasi
conformally flat and conformally symmetric 4-dimensional
\textit{LP}-Sasakian manifold to be proper biharmonic. Also, we
investigate proper biharmonic curves in the Lorentzian sphere
$S^{4}_{1}$.}

\\
\M{58C40, 53C42, 53C25.}

\section{Introduction}

The theory of biharmonic functions is an old and rich subject. Biharmonic
functions have been studied since 1862 by Maxwell and Airy to describe a
mathematical model of elasticity. The theory of polyharmonic functions was
developed later on, for example, by E. Almansi, T. Levi-Civita and M.
Nicolescu. Recently, biharmonic functions on Riemannian manifolds were
studied by R. Caddeo and L. Vanchke \cite{10,17}, L. Sario, M. Nakai and C.
Wang \cite{52}.

In the last decade there has been a growing interest in the theory
of biharmonic maps which can be divided in two main research
directions. On the one side, constructing the examples and
classification results have become important from the differential
geometric aspect. The other side is the analytic aspect from the
point of view of partial differential equations (see
\cite{18,36,55,56,57}), because biharmonic maps are solutions of a
fourth order strongly elliptic semilinear PDE.

Let $C^{\infty }(M,N)$ denote the space of smooth maps $\Psi
:(M,g)\rightarrow (N,h)$ \ between two Riemannian manifolds. A map
$\Psi \in C^{\infty }(M,N)$ is called \textit{harmonic }if it is a
critical point of the \textit{energy} functional
\begin{equation*}
E:C^{\infty }(M,N)\rightarrow R,E(\Psi )=\frac{1}{2}\int_{M}|d\Psi |^{2}v_{g}
\end{equation*}
and is characterized by the vanishing of the tension field $\tau
(\Psi )=trace\nabla d\Psi$ where $\nabla$ is a connection induced
from the Levi-Civita connection $\nabla^{M}$ of $M$ and the
pull-back connection $\nabla^{\Psi}$. As a generalization of
harmonic maps, biharmonic maps between Riemannian manifolds were
introduced by J. Eells and J. H. Sampson in
\cite{Eells-Sampson}. \textit{Biharmonic maps} between Riemannian manifolds $%
\Psi :(M,g)\rightarrow (N,h)$ are the critical points of the \textit{%
bienergy functional}
\begin{equation*}
E_{2}(\Psi )=\frac{1}{2}\int_{M}|\tau (\Psi )|^{2}v_{g}.
\end{equation*}
The first variation formula for the bienergy \ which is derived in \cite%
{Jiang1, Jiang2} shows that the Euler-Lagrange equation for the bienergy is
\begin{equation*}
\tau _{2}(\Psi )=-J(\tau (\Psi ))=-\Delta \tau (\Psi
)-traceR^{N}(d\Psi ,\tau (\Psi ))d\Psi =0,
\end{equation*}%
where $\Delta =-trace(\nabla ^{\Psi }\nabla ^{\Psi }-\nabla _{\nabla }^{\Psi
})$ is the rough Laplacian on the sections of $\Psi ^{-1}TN$ and $%
R^{N}(X,Y)=[\nabla _{X},\nabla _{Y}]-\nabla _{\lbrack X,Y]\text{ }}$ is the
curvature operator on $N$. From the expression of the bitension field $\tau
_{2}$, it is clear that a harmonic map is automatically a biharmonic map. So
non-harmonic biharmonic maps which are called proper biharmonic maps are
more interesting.

In a different setting, B. Y. Chen \cite{Chen} defined biharmonic
submanifolds $M\subset R^{n}$ of the Euclidean space as those with
harmonic mean curvature vector field, that is $\Delta H=0,$ where
$\Delta $ is the rough Laplacian, and stated the following

\begin{itemize}
\item Conjecture: Any biharmonic submanifold of the Euclidean
space is harmonic, that is minimal.
\end{itemize}

If the definition \ of biharmonic maps is applied to Riemannian
immersions into Euclidean space, the notion of Chen's biharmonic
submanifold is obtained, so the two definitions agree.

The non-existence theorems for the case of non-positive sectional curvature
codomains, as well as the
\begin{itemize}
\item Generalized Chen's conjecture: Biharmonic submanifolds of a manifold $%
N $ with $Riem^{N}\leq 0$ are minimal,
\end{itemize}
encouraged the study of proper biharmonic submanifolds, that is
submanifolds such that the inclusion map is a biharmonic map, in
spheres or another non-negatively curved spaces (see \cite{A2, A3,
A10, A11, A13, A16}).

Of course, the first and easiest examples can be found by looking at
differentiable curves in a Riemannian manifold. Obviously geodesics
are biharmonic. Non-geodesic biharmonic curves are called proper
biharmonic curves. Chen and Ishikawa \cite{D13} showed non-existence
of proper biharmonic curves in Euclidean 3-space $E^{3}.$ Moreover
they classified all proper biharmonic curves in Minkowski 3-space
$E_{1}^{3}$ (see also \cite{A12}). Caddeo, Montaldo and Piu showed
that on a surface with non-positive Gaussian curvature, any
biharmonic curve is a geodesic of the surface \cite{C14}. So they
gave a positive answer to generalized Chen's conjecture. Caddeo et
al. in \cite{A2} studied biharmonic curves in the unit 3-sphere.
More precisely, they showed that proper biharmonic curves in $S^{3}
$ are circles of geodesic curvature 1 or helices which are geodesics
in the Clifford minimal torus. Then the same authors studied the
biharmonic submanifolds of unit n-sphere \cite{A3}.

On the other hand, there are a few results on biharmonic curves in
arbitrary Riemannian manifolds. The biharmonic curves in the
Heisenberg group $H_{3}$ are investigated in \cite{D10} by Caddeo et
al. In \cite{A10} Fetcu studied biharmonic curves in the generalized
Heisenberg group and obtained two families of proper biharmonic
curves. Also, the explicit parametric equations for the biharmonic
curves on Berger spheres $S_{\varepsilon }^{3}$ are obtained by
Balmu\c{s} in \cite{C6}.

A generalization of Riemannian manifolds with constant sectional
curvature is represented by Sasakian space forms. In particular, a
simply connected three-dimensional Sasakian space form of constant
holomorphic sectional curvature $1$ is isometric to $S^{3}.$ So in
this context J. Inoguchi classified in \cite{A11} the proper
biharmonic Legendre curves and Hopf cylinders in a $3$-dimensional
Sasakian space form and in \cite{B11} the explict parametric
equations were obtained. T. Sasahara \cite{C53}, analyzed the proper
biharmonic Legendre surfaces in Sasakian space forms and in the case
when the ambient space is the unit $5$-dimensional sphere $S^{5}$ he
obtained their explicit representations.

Other results on biharmonic Legendre curves and biharmonic
anti-invariant surfaces in Sasakian space forms and $(\kappa ,\mu
)$-manifolds are given in \cite{C2,C1}.

In this paper we give some necessary and sufficient condition for a
spacelike and a timelike curve lying in a 4-dimensional conformally flat,
quasi conformally flat and conformally symmetric Lorentzian para-Sasakian
manifold to be proper biharmonic.

The study of Lorentzian almost paracontact manifolds was initiated
by Matsumoto in 1989 \cite{F11}. Also he introduced the notion of
Lorentzian para-Sasakian ( for short \textit{LP}-Sasakian )
manifold. I. Mihai and R. Rosca \cite{H3} defined the same notion
independently and thereafter many authors \cite{H4,H5,mukut} studied
\textit{LP}-Sasakian manifolds.

\section{Preliminaries}

\setcounter{equation}{0} \renewcommand{\theequation}{2.\arabic{equation}}

\subsection{Biharmonic maps between Riemannian manifolds}

\setcounter{equation}{0} \renewcommand{\theequation}{2.1.\arabic{equation}}
Let $(M,g)$ and $(N,h)$ be Riemannian manifolds and $\Psi :(M,g)\rightarrow
(N,h)$ \ be a smooth map. The tension field of $\Psi $ is given by $\tau
(\Psi )=trace\nabla d\Psi $, and for any compact domain $\Omega \subseteq M$%
, the bienergy is defined by
\begin{eqnarray*}
E_{2}(\Psi )=\frac{1}{2}\int_{\Omega }|\tau (\Psi )|^{2}v_{g}.
\end{eqnarray*}%
Then a smooth map $\Psi $ is called biharmonic map if it is a critical point
of the bienergy functional for any compact domain $\Omega \subseteq M.$ The
first variation formula for the bienergy functional is given by
\begin{eqnarray*}
\frac{dE_{2}(\Psi _{t})}{dt}|_{t=0}=\int_{\Omega }<\tau _{2}(\Psi ),w>v_{g},
\end{eqnarray*}%
where $v_{g}$ is the volume element, $w$ is the variational vector field
associated to the variation $\{\Psi _{t}\}$ of $\Psi $ and
\begin{eqnarray*}
\tau _{2}(\Psi )=-J(\tau _{2}(\Psi ))=-\Delta ^{\Psi }\tau (\Psi
)-traceR^{N}(d\Psi ,\tau (\Psi ))d\Psi .
\end{eqnarray*}%
Here $\Delta ^{\Psi }$ is the rough Laplacian on the sections of the
pull-back bundle $\Psi ^{-1}TN$ which is defined by
\begin{eqnarray*}
\Delta ^{\Psi }V=-\sum_{i=1}^{m}\{\nabla _{e_{i}}^{\Psi }\nabla
_{e_{i}}^{\Psi }V-\nabla _{\nabla _{e_{i}}^{M}e_{i}}^{\Psi
}V\},\,\,\,\,\,V\in \Gamma (\Psi ^{-1}TN),
\end{eqnarray*}%
where $\nabla $ is the pull-back connection on the pull-back bundle $\Psi
^{-1}TN$ and $\{e_{i}\}_{i=1}^{m}$ is an ortonormal frame on $M.$

From the definition of bienergy and the equation $\tau _{2}(\Psi )$, some
remarks on biharmonic maps are following:

\begin{itemize}
\item a map $\Psi $ is biharmonic if and only if its tension field is in the
kernel of the Jacobi operator;

\item a harmonic map is obviously a biharmonic map;

\item a harmonic map is an absolute minimum of the bienergy.
\end{itemize}

In particular, if the target manifold $N$ is the Euclidean space
$E^{m}$,
then the biharmonic equation of a map $\Psi :M\rightarrow E^{m}$ is%
\begin{equation*}
\Delta ^{2}\Psi =0,
\end{equation*}%
where $\Delta $ is the Laplace-Beltrami operator of $(M,g).$ Also,
biharmonic parametrized curves $\gamma:I\subset{R}\rightarrow M$ are
solutions of the fourth order differential equation
\begin{eqnarray*}
\nabla_{T}^{3}T-R(T,\nabla_{T}T)T=0.
\end{eqnarray*}

\subsection{Lorentzian Almost paracontact manifolds}

\setcounter{equation}{0} \renewcommand{\theequation}{2.2.\arabic{equation}}
Let $M$ \ be an $n$-dimensional smooth connected paracompact Hausdroff
manifold with a Lorentzian metric $g$, i.e., $g$ is a smooth symmetric
tensor field of type $(0,2)$ such that at every point $p\in M$, the tensor $%
g_{p}:T_{p}M\times T_{p}M\rightarrow R$ is a non-degenerate inner product of
signature $(-,+,...,+),$ where $T_{p}M$ is the tangent space of $M$ at the
point $p.$ Then $(M,g)$ is known to be a Lorentzian manifold. A non-zero
vector $X_{p}\in T_{p}M$ can be spacelike, null or timelike ,if it satisfies
$g_{p}(X_{p},X_{p})\geq 0,$ $g_{p}(X_{p},X_{p})=0$ $(X_{p}\neq 0)$ or $%
g_{p}(X_{p},X_{p})< 0$ respectively.

Let $M$ be an n-dimensional differentiable manifold equipped with a triple $%
(\phi ,\xi ,\eta ),$ where $\phi $ is a $(1,1)$ tensor field, $\xi $ is a
vector field, $\eta $ is a 1-form on $M$ such that \cite{F11}
\begin{eqnarray}
\eta (\xi ) &=&-1, \\
\phi ^{2} &=&I+\eta \otimes \xi ,
\end{eqnarray}%
where $I$ denotes the identity map of $T_{p}M$ and $\otimes $ is the tensor
product. The equations (2.2.1) and (2.2.2) imply that%
\begin{eqnarray*}
\eta \circ \phi &=&0, \\
\phi \xi &=&0, \\
rank(\phi ) &=&n-1.
\end{eqnarray*}%
Then $M$ admits a Lorentzian metric $g$, such that%
\begin{eqnarray*}
g(\phi X,\phi Y)=g(X,Y)+\eta (X)\eta (Y),
\end{eqnarray*}%
and $M$ is said to admit a Lorentzian almost paracontact structure
$(\phi ,\xi ,\eta ,g).$ Then we get
\begin{eqnarray*}
g(X,\xi ) &=&\eta (X), \\
\Phi (X,Y) &\equiv &g(X,\phi Y)\equiv g(\phi X,Y)\equiv \Phi (Y;X), \\
(\nabla _{X}\Phi )(Y,Z) &=&g(Y,(\nabla _{X}\phi )Z)=(\nabla _{X}\Phi )(Z,Y),
\end{eqnarray*}%
where $\nabla $ is the covariant differentiation with respect to $g$. It is
clear that Lorentzian metric $g$ makes $\xi $ a timelike unit vector field,
i.e, $g(\xi ,\xi )=-1.$ The manifold $M$ equipped with a Lorentzian almost
paracontact structure $(\phi ,\xi ,\eta ,g)$ is called a Lorentzian almost \
paracontact manifold (for short $LAP$-manifold) \cite{F11,F12}.

In equations (2.2.1) and (2.2.2) if we replace $\xi $ by $-\xi $, we obtain
an almost paracontact structure on $M$ defined by Sat\={o} \cite{F25}.

A Lorentzian almost \ paracontact manifold $M$ endowed with the structure $%
(\phi ,\xi ,\eta ,g)$ is called a Lorentzian paracontact manifold ( for
short \textit{LP}-manifold) \cite{F11} if%
\begin{eqnarray*}
\Phi (X,Y)=\frac{1}{2}((\nabla _{X}\eta )Y+(\nabla _{Y}\eta )X.
\end{eqnarray*}
\ \

A Lorentzian almost \ paracontact manifold $M$ endowed with the structure $%
(\phi ,\xi ,\eta ,g)$ is called a Lorentzian para Sasakian manifold ( for
short \textit{LP}-Sasakian) \cite{F11} if%
\begin{eqnarray*}
(\nabla _{X}\phi )Y=g(\phi X,\phi Y)\xi +\eta (Y)\phi ^{2}X,
\end{eqnarray*}%
or equivalently,%
\begin{eqnarray*}
(\nabla _{X}\phi )Y=\eta (Y)X+g(X;Y)\xi +2\eta (X)\eta (Y)\xi,
\end{eqnarray*}%
or equivalently,%
\begin{eqnarray*}
\, \, \, \, \, \, \, \, \, \, \, \,\, \, \, \, \, \,\,\, \, \, \, \,
\,(\nabla _{X}\Phi )(Y,Z)=g(X,Y)\eta (Z)+g(X;Z)\eta (Y)+2\eta (X)\eta
(Y)\eta (Z).
\end{eqnarray*}%
In a \textit{LP}-Sasakian manifold the $1$-form $\eta $ is closed.

Also Matsomoto in \cite{F11} showed that if an $n$-dimensional Lorentzian
manifold $(M,g)$ admits a timelike unit vector field $\xi $ such that the $1$%
-form $\eta $ associated to $\xi $ is closed and satisfies
\begin{eqnarray*}
\, \, \, \, \, \, \, \, \, \, \, \,\, \, \, \, \, \,\,\, \, \, \, \,
\,(\nabla _{X}\nabla _{Y}\eta )Z=g(X,Y)\eta (Z)+g(X,Z)\eta (Y)+2\eta (X)\eta
(Y)\eta (Z),
\end{eqnarray*}%
then $(M,g)$ admits an \textit{LP}-Sasakian structure.

An \textit{LP}-Sasakian manifold $M^{n}$ is said to be $\eta
$-Einstein if its Ricci tensor $S$ is of the form
\begin{eqnarray*}
S(X,Y)=ag(X,Y)+b\eta (X)\eta (Y),X,Y\in \Gamma (TM),
\end{eqnarray*}%
where $a$ and $b$ are functions on $M^{n}$ \cite{G10,G9}.

The conformal curvature tensor $C$ is defined by
\begin{eqnarray*}
C(X,Y)Z &=&R(X,Y)Z-\frac{1}{n-2}\{g(Y,Z)QX-g(X,Z)QY  \notag \\
&&+S(Y,Z)X-S(X,Z)Y\}+\frac{r}{(n-1)(n-2)}\{g(Y,Z)X-g(X,Z)Y\},
\end{eqnarray*}%
where $S(X,Y)=g(QX,Y).$ If $C=0$ then the \textit{LP}-Sasakian
manifold is called conformally flat.

The quasi-conformal curvature tensor $\tilde{C}$ is given by%
\begin{eqnarray*}
\widetilde{C}(X,Y)Z &=&aR(X,Y)Z+b\{S(Y,Z)X-S(X,Z)Y+g(Y,Z)QX  \notag \\
&&-g(X,Z)QY\}-\frac{r}{n}(\frac{a}{n-1}+2b)\{g(Y,Z)X-g(X,Z)Y\},
\end{eqnarray*}%
where $a,$ $b$ constants such that $ab\neq 0$ and $S(Y,Z)=g(QY,Z)$. If $%
\widetilde{C}=0$ then the \textit{LP}-Sasakian manifold is called
quasi conformally flat. In \cite{H} it was proved that a conformally
flat and a quasi conformally flat \textit{LP}-Sasakian manifold is
of constant curvature and the value of this constant is +1. Also the
same authors showed in \cite{H} that if in an \textit{LP}-Sasakian
manifold $M^{n}$ $(n>3)$ the relation $R(X,Y).C=0$ holds, then it is
locally isometric to a Lorentzian unit sphere.

For a conformally symmetric Riemannian manifold \cite{H1}, we have
$\nabla C=0.$ Hence for such a manifold $R(X,Y).C=0$ holds. Thus a
conformally symmetric \textit{LP}-Sasakian manifold $M^{n}$ $(n>3)$
is locally isometric to a Lorentzian unit sphere \cite{H}.

For a conformally flat, quasi conformally flat and conformally
symmetric \textit{LP}-Sasakian manifold $M^{n}$, we have \cite{H}
\begin{eqnarray}
\text{ \ \ \ \ \ \ \ \ \ \ \ \ }R(X,Y)Z=g(Y,Z)X-g(X,Z)Y,\text{ \ \ \ \ \ \ \
\ }X,\text{ }Y,\text{ }Z\in \Gamma (TM).
\end{eqnarray}

An arbitrary curve $\gamma :I\rightarrow M,$ $\gamma =\gamma (s),$
in a \textit{LP}-Sasakian manifold is called spacelike, timelike or
null (lightlike), if all of its velocity vectors $\gamma ^{\prime
}(s)$ are respectively spacelike, timelike or null (lightlike). If
$\gamma (s)$ is a spacelike or timelike curve, we can reparametrize
it such that $g(\gamma'(s),\gamma'(s))=\varepsilon$ where
$\varepsilon=1$ if $\gamma$ is spacelike and $\varepsilon=-1$ if
$\gamma$ is timelike, respectively. In this case $\gamma (s)$ is
said to be unit speed or arclenght parametrization.

Denote by $\{T(s),N(s),B_{1}(s),B_{2}(s)\}$ the moving Frenet frame
along the curve $\gamma (s)$ in a \textit{LP}-Sasakian manifold.
Then $T,N,B_{1},B_{2}$ are respectively, the tangent, the principal
normal, the first binormal and the second binormal vector fields. A
spacelike or timelike curve $\gamma (s)$ is said to be parametrized
by arclenght function $s$, if $g(\gamma ^{\prime }(s),\gamma
^{\prime }(s))=\pm 1.$

Let $\gamma (s)$ be a curve in \textit{LP}-Sasakian manifold
parametrized by arclenght function $s.$ Then for the curve $\gamma $
the following Frenet equations are given in \cite{il5}:

\textbf{Case I. }$\gamma $\textbf{\ is a spacelike curve:}

Then $T$ is a spacelike vector, so depending on the casual character of the
principal normal vector $N$ and the first binormal vector $B_{1}$, we have
the following Frenet formulas:\newline
\textbf{Case I.1.} $N$ and $B_{1}$ are spacelike;%
\begin{equation}
\label{1}
\left[
\begin{array}{c}
\nabla _{T}T \\
\nabla _{T}N \\
\nabla _{T}B_{1} \\
\nabla _{T}B_{2}%
\end{array}%
\right] =\left[
\begin{array}{cccc}
0 & k_{1} & 0 & 0 \\
-k_{1} & 0 & k_{2} & 0 \\
0 & -k_{2} & 0 & k_{3} \\
0 & 0 & k_{3} & 0%
\end{array}%
\right] \left[
\begin{array}{c}
T \\
N \\
B_{1} \\
B_{2}%
\end{array}%
\right] ,
\end{equation}%
where $T,$ $N,$ $B_{1},$ $B_{2}$ are mutually orthogonal vectors
satisfying (\ref{1})
the equations%
\begin{equation*}
g(T,T)=g(N,N)=g(B_{1},B_{1})=1,\text{ \ \ \ \ \ \ \ }g(B_{2},B_{2})=-1.
\end{equation*}%
\textbf{Case I.2.} $N$ is spacelike, $B_{1}$ is timelike;%
\begin{eqnarray}
\left[
\begin{array}{c}
\nabla _{T}T \\
\nabla _{T}N \\
\nabla _{T}B_{1} \\
\nabla _{T}B_{2}%
\end{array}%
\right] =\left[
\begin{array}{cccc}
0 & k_{1} & 0 & 0 \\
-k_{1} & 0 & k_{2} & 0 \\
0 & k_{2} & 0 & k_{3} \\
0 & 0 & k_{3} & 0%
\end{array}%
\right] \left[
\begin{array}{c}
T \\
N \\
B_{1} \\
B_{2}%
\end{array}%
\right] ,
\end{eqnarray}%
where $T,$ $N,$ $B_{1},$ $B_{2}$ are mutually orthogonal vectors
satisfying
the equations%
\begin{equation*}
g(T,T)=g(N,N)=g(B_{2},B_{2})=1,\text{ \ \ \ \ \ \ \ }g(B_{1},B_{1})=-1.
\end{equation*}%
\textbf{Case I.3.} $N$ is spacelike, $B_{1}$ is null;%
\begin{eqnarray}
\left[
\begin{array}{c}
\nabla _{T}T \\
\nabla _{T}N \\
\nabla _{T}B_{1} \\
\nabla _{T}B_{2}%
\end{array}%
\right] =\left[
\begin{array}{cccc}
0 & k_{1} & 0 & 0 \\
-k_{1} & 0 & k_{2} & 0 \\
0 & 0 & k_{3} & 0 \\
0 & -k_{2} & 0 & -k_{3}%
\end{array}%
\right] \left[
\begin{array}{c}
T \\
N \\
B_{1} \\
B_{2}%
\end{array}%
\right] ,
\end{eqnarray}%
where $T,$ $N,$ $B_{1},$ $B_{2}$ satisfy the equations%
\begin{eqnarray*}
g(T,T) &=&\text{\ }g(N,N)=1,\text{ \ \ \ \ \ \ \ }%
g(B_{1},B_{1})=g(B_{2},B_{2})=0, \\
g(T,N) &=&g(T,B_{1})=g(T,B_{2})=g(N,B_{1})=g(N,B_{2})=0,g(B_{1},B_{2})=1.
\end{eqnarray*}%
\textbf{Case I.4.} $N$ is timelike, $B_{1}$ is spacelike;%
\begin{eqnarray}
\left[
\begin{array}{c}
\nabla _{T}T \\
\nabla _{T}N \\
\nabla _{T}B_{1} \\
\nabla _{T}B_{2}%
\end{array}%
\right] =\left[
\begin{array}{cccc}
0 & k_{1} & 0 & 0 \\
k_{1} & 0 & k_{2} & 0 \\
0 & k_{2} & 0 & k_{3} \\
0 & 0 & -k_{3} & 0%
\end{array}%
\right] \left[
\begin{array}{c}
T \\
N \\
B_{1} \\
B_{2}%
\end{array}%
\right] ,
\end{eqnarray}%
where $T,$ $N,$ $B_{1},$ $B_{2}$ are mutually orthogonal vectors
satisfying
the equations%
\begin{equation*}
g(T,T)=g(B_{1},B_{1})=g(B_{2},B_{2})=1,\text{ \ \ \ \ \ \ \ \ \
}g(N,N)=-1.
\end{equation*}%
\textbf{Case I.5.} $N$ is null, $B_{1}$ is spacelike;%
\begin{eqnarray}
\left[
\begin{array}{c}
\nabla _{T}T \\
\nabla _{T}N \\
\nabla _{T}B_{1} \\
\nabla _{T}B_{2}%
\end{array}%
\right] =\left[
\begin{array}{cccc}
0 & k_{1} & 0 & 0 \\
0 & 0 & k_{2} & 0 \\
0 & k_{3} & 0 & -k_{2} \\
-k_{1} & 0 & -k_{3} & 0%
\end{array}%
\right] \left[
\begin{array}{c}
T \\
N \\
B_{1} \\
B_{2}%
\end{array}%
\right] ,
\end{eqnarray}%
where $T,$ $N,$ $B_{1},$ $B_{2}$ satisfy the equations%
$$
g(T,T) =\text{\
}g(B_{1},B_{1})=1,\,\,\,\,\,\,\,\,g(N,N)=g(B_{2},B_{2})=0,
$$
$$
g(T,N)
=g(T,B_{1})=g(T,B_{2})=g(N,B_{1})=g(B_{1},B_{2})=0,g(N,B_{2})=1.
$$

\textbf{Case II. }$\gamma $\textbf{\ is a timelike curve:}

In this case $T$ is a timelike vector, so the Frenet formulae have
the form
\textbf{\ \ \ }%
\begin{eqnarray}
\left[
\begin{array}{c}
\nabla _{T}T \\
\nabla _{T}N \\
\nabla _{T}B_{1} \\
\nabla _{T}B_{2}%
\end{array}%
\right] =\left[
\begin{array}{cccc}
0 & k_{1} & 0 & 0 \\
k_{1} & 0 & k_{2} & 0 \\
0 & -k_{2} & 0 & k_{3} \\
0 & 0 & -k_{3} & 0%
\end{array}%
\right] \left[
\begin{array}{c}
T \\
N \\
B_{1} \\
B_{2}%
\end{array}%
\right] ,
\end{eqnarray}%
where $T,$ $N,$ $B_{1},$ $B_{2}$ are mutually orthogonal vectors satisfying
the equations%
\begin{equation*}
g(N,N)=g(B_{1},B_{1})=g(B_{2},B_{2})=1,\text{ \ \ \ \ \ \ \
}g(T,T)=-1.
\end{equation*}

\section{Biharmonic curves in \textit{LP}-Sasakian manifolds}

\setcounter{equation}{0}
\renewcommand{\theequation}{3.\arabic{equation}} In this section we
characterize the spacelike and timelike proper biharmonic curves in
a 4-dimensional conformally flat, quasi conformally flat and
conformally symmetric Lorentzian para-Sasakian
(\textit{LP}-Sasakian) manifold.\\

\begin{theorem} \textit{Let $M$ be a 4-dimensional
conformally flat, quasi conformally flat or conformally symmetric
\textit{LP}-Sasakian manifold and $\gamma :I\rightarrow M$ be a
spacelike curve parametrized by arclength. Suppose that
$\{T,N,B_{1},B_{2}\}$
be an orthonormal Frenet frame field tangent to $M$ along $\gamma $ such that $%
g(T,T)=g(N,N)=g(B_{1},B_{1})=1$ and $g(B_{2},B_{2})=-1.$ Then $\gamma
:I\rightarrow M$ is a proper biharmonic curve if and only if either $\gamma $
is a circle with $k_{1}=1,$ or $\gamma $ is a helix with $%
k_{1}^{2}+k_{2}^{2}=1.$}\\
\end{theorem}
\begin{proof} Let $M$ be a $4$-dimensional conformally flat, quasi
conformally flat or conformally symmetric \textit{LP}-Sasakian
manifold endowed with the structure$\ (\phi ,\xi ,\eta ,g)$ and
$\gamma :I\rightarrow M$ be a curve parametrized by arclength.
Suppose that $\gamma $ is a spacelike curve that is its
velocity vector $T=\gamma ^{\prime }(s)$ is spacelike. Let $%
\{T,N,B_{1},B_{2}\}$ be an orthonormal Frenet frame field tangent to $M$ along $%
\gamma $, where $N$ is the unit spacelike vector field in the direction $%
\nabla _{T}T$ , $B_{1}$ is a unit spacelike and $B_{2}$ is a unit
timelike vector. The tension field of $\gamma $ is $\tau (\gamma
)=\nabla _{T}T.$ Then by using the Frenet formulas (2.2.4) and the
equation (2.2.3) we
obtain the Euler-Lagrange equation of the bienergy:%
\begin{eqnarray*}
\tau _{2}(\gamma ) &=&\nabla _{T}^{3}T-R(T,\nabla _{T}T)T \\
&=&\nabla _{T}^{3}T-R(T,k_{1}N)T \\
&=&(-3k_{1}k_{1}^{\prime })T+(k_{1}^{\prime \prime
}-k_{1}^{3}-k_{1}k_{2}^{2})N \\
&&+(2k_{1}^{\prime }k_{2}+k_{1}k_{2}^{\prime
})B_{1}+(k_{1}k_{2}k_{3})B_{2}-k_{1}R(T,N)T \\
&=&(-3k_{1}k_{1}^{\prime })T+(k_{1}^{\prime \prime
}-k_{1}^{3}-k_{1}k_{2}^{2}+k_{1})N \\
&&+(2k_{1}^{\prime }k_{2}+k_{1}k_{2}^{\prime })B_{1}+(k_{1}k_{2}k_{3})B_{2}
\\
&=&0.
\end{eqnarray*}
where $k_{1},$ $k_{2}$ and $k_{3}$ are respectively the first, the
second and the third curvature of the curve $\gamma (s).$

It follows that $\gamma $ is a biharmonic curve if and only if%
\begin{eqnarray*}
k_{1}k_{1}^{\prime } &=&0, \\
k_{1}^{\prime \prime }-k_{1}(k_{1}^{2}+k_{2}^{2}-1) &=&0 ,\\
2k_{1}^{\prime }k_{2}+k_{1}k_{2}^{\prime } &=&0 ,\\
k_{1}k_{2}k_{3} &=&0.
\end{eqnarray*}%
If we look for nongeodesic solutions , that is for biharmonic curves with $%
k_{1}\neq 0,$ we obtain%
\begin{eqnarray*}
k_{1}&=&constant\neq 0,\,k_{2}=constant, \\
k_{1}^{2}+k_{2}^{2} &=&1, \\
k_{2}k_{3} &=&0.
\end{eqnarray*}%
This completes the proof.
\end{proof}

\begin{theorem}
Let $M$ be a 4-dimensional conformally flat, quasi conformally flat
or conformally symmetric \textit{LP}-Sasakian manifold and $\gamma
:I\rightarrow M$ be a spacelike curve parametrized by arclength.
Suppose that $\{T,N,B_{1},B_{2}\}$
be an orthonormal Frenet frame field tangent to $M$ along $\gamma $ such that $%
g(T,T)=g(N,N)=g(B_{2},B_{2})=1$ and $g(B_{1},B_{1})=-1.$ Then $\gamma
:I\rightarrow M$ is a proper biharmonic curve if and only if either $\gamma $
is a circle with $k_{1}=1,$ or $\gamma $ is a helix with $%
k_{1}^{2}-k_{2}^{2}=1.$
\end{theorem}

\begin{proof}
Let $M$ be a $4$-dimensional conformally flat, quasi conformally
flat or conformally symmetric \textit{LP}-Sasakian manifold endowed
with the structure$\ (\phi ,\xi ,\eta ,g)$ and $\gamma :I\rightarrow
M$ be a curve parametrized by arclength. Suppose that $\gamma $ is a
spacelike curve that is its
velocity vector $T=\gamma ^{\prime }(s)$ is spacelike. Let $%
\{T,N,B_{1},B_{2}\}$ be an orthonormal Frenet frame field tangent to $M$ along $%
\gamma $, where $N$ is the unit spacelike vector field in the direction $%
\nabla _{T}T$ , $B_{2}$ is a unit spacelike and $B_{1}$ is a unit
timelike vector. Since the tension field of $\gamma $ is $\tau
(\gamma )=\nabla _{T}T$ then by using the Frenet formulas given in
(2.2.5) and the equation (2.2.3), we obtain the biharmonic equation
for $\gamma$:
\begin{eqnarray*}
\tau _{2}(\gamma ) &=&\nabla _{T}^{3}T-R(T,\nabla _{T}T)T \\
&=&\nabla _{T}^{3}T-R(T,k_{1}N)T \\
&=&(-3k_{1}k_{1}^{\prime })T+(k_{1}^{\prime \prime
}-k_{1}^{3}+k_{1}k_{2}^{2})N \\
&&+(2k_{1}^{\prime }k_{2}+k_{1}k_{2}^{\prime
})B_{1}+(k_{1}k_{2}k_{3})B_{2}-k_{1}R(T,N)T \\
&=&(-3k_{1}k_{1}^{\prime })T+(k_{1}^{\prime \prime
}-k_{1}^{3}+k_{1}k_{2}^{2}+k_{1})N \\
&&+(2k_{1}^{\prime }k_{2}+k_{1}k_{2}^{\prime })B_{1}+(k_{1}k_{2}k_{3})B_{2}
\\
&=&0.
\end{eqnarray*}
where $k_{1},$ $k_{2}$ and $k_{3}$ are respectively the first, the second
and the third curvature of curve $\gamma (s).$

It follows that $\gamma $ is a biharmonic curve if and only if%
\begin{eqnarray*}
k_{1}k_{1}^{\prime } &=&0, \\
k_{1}^{\prime \prime }-k_{1}(k_{1}^{2}-k_{2}^{2}-1) &=&0 ,\\
2k_{1}^{\prime }k_{2}+k_{1}k_{2}^{\prime } &=&0, \\
k_{1}k_{2}k_{3} &=&0.
\end{eqnarray*}%
If we look for nongeodesic solutions , that is for biharmonic curves with $%
k_{1}\neq 0,$ we obtain%
\begin{eqnarray*}
k_{1} &=&cons\tan t\neq 0,k_{2}=cons\tan t ,\\
k_{1}^{2}-k_{2}^{2} &=&1, \\
k_{2}k_{3} &=&0.
\end{eqnarray*}%
This completes the proof.
\end{proof}

\begin{theorem}
Let $M$ be a 4-dimensional conformally flat, quasi conformally flat
or conformally symmetric \textit{LP}-Sasakian manifold and $\gamma
:I\rightarrow M$ be a spacelike curve parametrized by arclength.
Suppose that $\{T,N,B_{1},B_{2}\}$ be a moving Frenet frame such
that $N$ is a spacelike and $B_{1}$ is a null vector. Then $\gamma
:I\rightarrow M$ is a proper biharmonic curve if and only if $k_{1}=
1$ and $\ln k_{2}(s)=-\int k_{3}(s)\,ds.$
\end{theorem}

\begin{proof}
Let $\gamma :I\rightarrow M$ be a spacelike curve parametrized by
arclength on a 4-dimensional conformally flat, quasi conformally
flat or conformally symmetric \textit{LP}-Sasakian manifold $M.$
Suppose that $\{T,N,B_{1},B_{2}\}$ be a moving Frenet frame such
that
$$
g(T,T) =\text{\ }g(N,N)=1,\text{ \ \ \ \ \ \ \ }%
g(B_{1},B_{1})=g(B_{2},B_{2})=0,
$$
$$
g(T,N) =g(T,B_{1})=g(T,B_{2})=g(N,B_{1})=g(N,B_{2})=0,\text{ \ \ \ \ \ \ }%
g(B_{1},B_{2})=1.
$$

Then by using the Frenet equations given by (2.2.6), we have%
\begin{eqnarray*}
\tau _{2}(\gamma ) &=&\nabla _{T}^{3}T-R(T,\nabla _{T}T)T \\
&=&\nabla _{T}^{3}T-R(T,k_{1}N)T \\
&=&(-3k_{1}k_{1}^{\prime })T+(k_{1}^{\prime \prime }-k_{1}^{3}+k_{1})N \\
&&+(2k_{1}^{\prime }k_{2}+k_{1}k_{2}^{\prime }+k_{1}k_{2}k_{3})B_{1}
\end{eqnarray*}%
where $k_{1},$ $k_{2}$ and $k_{3}$ are respectively the first, the second
and the third curvature of curve $\gamma (s).$ From the biharmonic equation
of $\gamma $ above, we can say $\gamma $ is a biharmonic curve if and only
if
\begin{eqnarray*}
k_{1}k_{1}^{\prime } &=&0, \\
k_{1}^{\prime \prime }-k_{1}^{3}+k_{1} &=&0, \\
2k_{1}^{\prime }k_{2}+k_{1}k_{2}^{\prime }+k_{1}k_{2}k_{3} &=&0.
\end{eqnarray*}%
For biharmonic curves with $k_{1}\neq 0$ that is if we investigate the
nongeodesic solutions, we obtain%
\begin{eqnarray*}
k_{1} &=&\mp 1, \\
k_{2}^{\prime }+k_{2}k_{3} &=&0.
\end{eqnarray*}
Thus we have the assertion of the theorem.
\end{proof}

\begin{theorem}
Let $M$ be a 4-dimensional conformally flat, quasi conformally flat
or conformally symmetric \textit{LP}-Sasakian manifold and $\gamma
:I\rightarrow M$ be a spacelike curve parametrized by arclength.
Suppose that $\{T,N,B_{1},B_{2}\}$
be an orthonormal Frenet frame field tangent to $M$ along $\gamma $ such that $%
g(T,T)=g(B_{1},B_{1})=g(B_{2},B_{2})=1$ and $g(N,N)=-1.$ Then $\gamma
:I\rightarrow M$ is a biharmonic curve if and only it is a geodesic of $M.$
\end{theorem}

\begin{proof}
Let $M$ be a $4$-dimensional conformally flat, quasi conformally
flat or conformally symmetric \textit{LP}-Sasakian manifold endowed
with the structure$\ (\phi ,\xi ,\eta ,g)$ and $\gamma :I\rightarrow
M$ be a curve parametrized by arclength. Suppose that $\gamma $ is a
spacelike curve that is its
velocity vector $T=\gamma ^{\prime }(s)$ is spacelike. Let $%
\{T,N,B_{1},B_{2}\}$ be an orthonormal Frenet frame field tangent to $M$ along $%
\gamma $, where $N$ is the unit timelike vector field in the direction $%
\nabla _{T}T$ , $B_{1}$ and $B_{2}$ are unit spacelike vectors. The
tension field of $\gamma $ is $\tau (\gamma )=\nabla _{T}T.$ Then by
using the tension field of $\gamma$, Frenet formulas in (2.2.7) and
the equation (2.2.3) we obtain the Euler-Lagrange equation of the
bienergy:%
\begin{eqnarray*}
\tau _{2}(\gamma ) &=&\nabla _{T}^{3}T-R(T,\nabla _{T}T)T \\
&=&\nabla _{T}^{3}T-R(T,k_{1}N)T \\
&=&(3k_{1}k_{1}^{\prime })T+(k_{1}^{\prime \prime
}+k_{1}^{3}+k_{1}k_{2}^{2})N \\
&&+(2k_{1}^{\prime }k_{2}+k_{1}k_{2}^{\prime
})B_{1}+(k_{1}k_{2}k_{3})B_{2}-k_{1}R(T,N)T \\
&=&(3k_{1}k_{1}^{\prime })T+(k_{1}^{\prime \prime
}+k_{1}^{3}+k_{1}k_{2}^{2}+k_{1})N\\
&&+(2k_{1}^{\prime }k_{2}+k_{1}k_{2}^{\prime
})B_{1}+(k_{1}k_{2}k_{3})B_{2} \\
&=&0.
\end{eqnarray*}

It follows that $\gamma $ is a biharmonic curve if and only if%
\begin{eqnarray*}
k_{1}k_{1}^{\prime } &=&0, \\
k_{1}^{\prime \prime }+k_{1}(k_{1}^{2}+k_{2}^{2}+1) &=&0, \\
2k_{1}^{\prime }k_{2}+k_{1}k_{2}^{\prime } &=&0 ,\\
k_{1}k_{2}k_{3} &=&0.
\end{eqnarray*}%
If we look for nongeodesic solutions , that is for biharmonic curves with $%
k_{1}\neq 0,$ we obtain%
\begin{eqnarray*}
k_{1} &=&cons\tan t\neq 0,k_{2}=cons\tan t, \\
k_{1}^{2}+k_{2}^{2} &=&-1, \\
k_{2}k_{3} &=&0.
\end{eqnarray*}%
This shows that we have no nongeodesic solution for the biharmonic equation
for the curve $\gamma .$
\end{proof}

\begin{theorem}
Let $M$ be a 4-dimensional conformally flat, quasi conformally flat
or conformally symmetric \textit{LP}-Sasakian manifold and $\gamma
:I\rightarrow M$ be a spacelike curve parametrized by arclength.
Suppose that $\{T,N,B_{1},B_{2}\}$ be a moving Frenet frame along
$\gamma$ such that $N$ is a null vector. Then $\gamma :I\rightarrow
M$ is a biharmonic curve if and only if $\gamma$ is a geodesic of
$M$.
\end{theorem}
\begin{proof}
Let $\gamma :I\rightarrow M$ be a spacelike curve parametrized by
arclength on a 4-dimensional conformally flat, quasi conformally
flat or conformally symmetric \textit{LP}-Sasakian manifold $M.$
Suppose that $\{T,N,B_{1},B_{2}\}$ be a moving Frenet frame along
the curve $\gamma$ such that
$$ g(T,T) =\text{\
}g(B_{1},B_{1})=1,\text{\ \ \ \ \ }g(N,N)=g(B_{2},B_{2})=0,
$$
$$
g(T,N) =g(T,B_{1})=g(T,B_{2})=g(N,B_{1})=g(B_{1},B_{2})=0,\text{\ \
\ \ \ }g(N,B_{2})=1.
$$
If we consider the Frenet formulas given in (2.2.8), we obtain the
biharmonic equation for the curve $\gamma$:
\begin{eqnarray*}
0=\tau _{2}(\gamma )&=&(k_{1}''+k_{1}k_{2}k_{3}+k_{1})N\\
&&+(2k_{1}'k_{2}+k_{1}k_{2}')B_{1}+(-k_{1}k_{2}^{2})B_{2}
\end{eqnarray*}
Then $\gamma$ is a biharmonic curve if and only if
\begin{eqnarray*}
k_{1}''+k_{1}k_{2}k_{3}+k_{1}&=&0,\\
2k_{1}'k_{2}+k_{1}k_{2}'&=&0,\\
k_{1}k_{2}^{2}&=&0.
\end{eqnarray*}
Since $\gamma$ is a spacelike curve with a null normal vector,
$k_{1}$ can take only two values: $0$ and $1$. If we look for
nongeodesic solutions, we get $k_{2}=0$. But from the first equation
above, we have a contradiction such that $k_{2}k_{3}+1=0$. So the
only biharmonic spacelike curves on $M$ with a null normal vector
are the geodesics of $M$.
\end{proof}

Let $M$ be a 4-dimensional conformally flat, quasi conformally flat
or conformally symmetric \textit{LP}-Sasakian manifold. Since $M$ is
locally isometric to a Lorentzian unit sphere $S_{1}^{4}$, by using
the above theorems we shall give some characterizations for
nongeodesic biharmonic curves in $S_{1}^{4}$:
\begin{proposition}
Let $\gamma :I\rightarrow S_{1}^{4}$ be a spacelike nongeodesic
biharmonic curve parametrized by arclenght and $\{T,N,B_{1},B_{2}\}$
be a Frenet frame along $\gamma$ such that the principal normal
vector N and first binormal vector $B_{1}$ are spacelike. Then
\begin{eqnarray}
\gamma^{(IV)}+2\gamma''+(1-k_{1}^{2})\gamma=0.
\end{eqnarray}
\end{proposition}
\begin{proof}
From the Frenet formulas (2.2.4), by taking the covariant derivative
of $\nabla_{T}N$ with respect to $T$, we have
\begin{eqnarray*}
\nabla^{2}_{T}N&=&-k_{1}\nabla_{T}T+k_{2}\nabla_{T}B_{1}\\
&&=-k_{1}^{2}N+k_{2}(-k_{2}N+k_{3}B_{2})\\
&&=-(k_{1}^{2}+k_{2}^{2})N+k_{2}k_{3}B_{2}\\
&&=-N.
\end{eqnarray*}
If we use the Gauss equation of $S_{1}^{4}\subset R_{1}^{5}$, that
for any vector field $X$ along $\gamma$ is
\begin{eqnarray*}
\nabla_{T}X=X'+<T,X>\gamma,
\end{eqnarray*}
we get
\begin{eqnarray*}
\nabla^{2}_{T}N&=&\nabla_{T}[N'+<T,N>\gamma] \\
&&=\nabla_{T}N'\\
&&=N''+<T,N'>\gamma\\
&&=N''+<T,\nabla_{T}N-<N,T>\gamma>\gamma\\
&&=N''+<T,\nabla_{T}N>\gamma\\
&&=N''-k_{1}\gamma
\end{eqnarray*}
and
\begin{eqnarray*}
N=\frac{1}{k_{1}}(\gamma''+\gamma).
\end{eqnarray*}
By substituting the above expressions of $\nabla^{2}_{T}N$ and $N$
in the equation $\nabla^{2}_{T}N+N=0$, we obtain the differential
equation (3.1).
\end{proof}

From Proposition 3.6, it is obvious that to find nongeodesic
biharmonic curves in $S_{1}^{4}$ we must investigate the solutions
of (3.1). By integrating the differential equation (3.1), we have

\begin{theorem}
Let $\gamma :I\rightarrow S_{1}^{4}$ be a spacelike nongeodesic
biharmonic curve parametrized by arclenght and $\{T,N,B_{1},B_{2}\}$
be a Frenet frame along $\gamma$ such that the principal normal
vector N and first binormal vector $B_{1}$ are spacelike. Then we
have two cases:
\begin{itemize}
\item $\gamma$ is a circle of radius $\frac{1}{\sqrt{2}}$;
\item $\gamma(s)=(0,\frac{\cos (as)}{\sqrt{2}},\frac{\sin
(as)}{\sqrt{2}},\frac{\cos (bs)}{\sqrt{2}},\frac{\sin
(bs)}{\sqrt{2}})$.
\end{itemize}
\end{theorem}

\begin{proof} If $k_{1}= 1$, then the general solution of (3.1)
is
\begin{eqnarray*}
\gamma(s)=c_{1}+c_{2}s+c_{3}\cos(\sqrt{2}s)+c_{4}\sin(\sqrt{2}s).
\end{eqnarray*}
Since $|\gamma|^{2}=1$ and $|\gamma'|^{2}=1$, we have $c_{2}=0$,
while $c_{1},\,c_{3},\,c_{4}$ are constant vectors ortogonal to each
other with $|c_{1}|^{2}=|c_{3}|^{2}=|c_{4}|^{2}=\frac{1}{2}$. Then
the solution becomes
\begin{eqnarray*}
\gamma(s)=(d_{1},\frac{\cos(\sqrt{2}s)}{\sqrt{2}},\frac{\sin(\sqrt{2}s)}{\sqrt{2}},d_{2},d_{3}),
\end{eqnarray*}
with $-d_{1}^{2}+d_{2}^{2}+d_{3}^{2}=\frac{1}{2}$. It is obvious
that $\gamma$ is a circle of radius $\frac{1}{\sqrt{2}}$.\\
If $0<k_{1}<1$, then the general solution of (3.1) is
\begin{eqnarray*}
\gamma(s)=c_{1}\cos(as)+c_{2}\sin(as)+c_{3}\cos(bs)+c_{4}\sin(bs)
\end{eqnarray*}
where $a=\sqrt{1-k_{1}}$ and $b=\sqrt{1+k_{1}}$. Since
$|\gamma|^{2}=1$ and $|\gamma'|^{2}=1$, we obtain that the vectors
$c_{i}$, $i=1,\,2,\,3,\,4$, are orthogonal to each other and
$|c_{1}|^{2}=|c_{1}|^{2}=|c_{3}|^{2}=|c_{4}|^{2}=\frac{1}{2}$. Then
the curve $\gamma$ becomes
\begin{eqnarray*}
\gamma(s)=(0,\frac{\cos (as)}{\sqrt{2}},\frac{\sin
(as)}{\sqrt{2}},\frac{\cos (bs)}{\sqrt{2}},\frac{\sin
(bs)}{\sqrt{2}}).
\end{eqnarray*}
\end{proof}

\begin{proposition}
Let $\gamma :I\rightarrow S_{1}^{4}$ be a spacelike nongeodesic
biharmonic curve parametrized by arclenght and $\{T,N,B_{1},B_{2}\}$
be a Frenet frame along $\gamma$ such that the principal normal
vector N is spacelike and first binormal vector $B_{1}$ is timelike.
Then
\begin{eqnarray}
\gamma^{(IV)}+2\gamma''+(1-k_{1}^{2})\gamma=0.
\end{eqnarray}
\end{proposition}

If $k_{1}= 1$, then it is obvious that the general solution of (3.2)
is a circle of radius $\frac{1}{\sqrt{2}}$. If $k_{1}>1$, then the
general solution of ( 3.2) is
\begin{eqnarray*}
\gamma(s)=c_{1}e^{as}+c_{2}e^{-as}+c_{3}\cos(bs)+c_{4}\sin(bs)
\end{eqnarray*} with $a=\sqrt{k_{1}-1}$ and $b=\sqrt{k_{1}+1}$. Here $c_{i}$, $i=1,\,2,\,3,\,4$, are constant
vectors. Since $|\gamma|^{2}=1$ and $|\gamma'|^{2}=1$, by choosing
\begin{eqnarray*}
c_{1}&=&(1,0,0,0,1),\,\,c_{2}=(-1,\frac{\sqrt{7}}{4},0,0,-\frac{3}{4}),\\
c_{3}&=&(0,0,\frac{1}{2},\frac{1}{2},0),\,\,
c_{4}=(-\frac{\sqrt{7}}{\sqrt{2}},\frac{1}{\sqrt{2}},0,0,-\frac{\sqrt{7}}{\sqrt{2}}),\,\,
\end{eqnarray*}
such that
\begin{eqnarray*}
<c_{1},c_{1}>&=&<c_{2},c_{2}>=0,\\
<c_{3},c_{3}>&=&<c_{4},c_{4}>=\frac{3}{b^{2}},\\
<c_{1},c_{2}>&=&\frac{1}{a^{2}},\\
<c_{1},c_{3}>&=&<c_{1},c_{4}>=0,\\
<c_{2},c_{3}>&=&<c_{2},c_{4}>=0,\\
<c_{3},c_{4}>&=&0,
\end{eqnarray*} with $a=2,\,\,b=\sqrt{6}$, we obtain following special
solution of differential equation (3.2)
\begin{eqnarray*}
\gamma(s)&=&(e^{2s}-e^{-2s}-\frac{\sqrt{7}}{\sqrt{2}}\sin(\sqrt{6}s),\frac{\sqrt{7}}{4}e^{-2s}+\frac{1}{\sqrt{2}}\sin(\sqrt{6}s),\nonumber\\
&&\frac{1}{2}\cos(\sqrt{6}s),\frac{1}{2}\cos(\sqrt{6}s),e^{2s}-\frac{3}{4}e^{-2s}-\frac{\sqrt{7}}{\sqrt{2}}\sin(\sqrt{6}s)),
\end{eqnarray*}
which is a helix with $k_{1}=5$ and $k_{2}=2\sqrt{6}$.

\begin{proposition}

Let $\gamma :I\rightarrow S_{1}^{4}$ be a spacelike nongeodesic
biharmonic curve parametrized by arclenght and $\{T,N,B_{1},B_{2}\}$
be a moving Frenet frame along $\gamma$ such that the principal
normal vector N is spacelike and first binormal vector $B_{1}$ is
null. Then
\begin{eqnarray}
\gamma^{(IV)}+2\gamma''=0.
\end{eqnarray}
\end{proposition}

It can be easily seen that the general solution of differential
equation (3.3) is a circle of radius $\frac{1}{\sqrt{2}}$.

Now let us investigate the biharmonicity of a timelike curve in a
4-dimensional conformally flat, quasi conformally flat and
conformally symmetric Lorentzian para-Sasakian
(\textit{LP}-Sasakian) manifold. We have,

\begin{theorem}
Let $M$ be a 4-dimensional conformally flat, quasi conformally flat
or conformally symmetric \textit{LP}-Sasakian manifold and $\gamma
:I\rightarrow M$ be a timelike curve parametrized by arclength. Then
$\gamma :I\rightarrow M$ is a proper biharmonic curve if and only if
either $\gamma $ is a circle with $k_{1}=1,$ or $\gamma $ is a helix
with $k_{1}^{2}-k_{2}^{2}=1.$
\end{theorem}

\begin{proof}
Let $M$ be a $4$-dimensional conformally flat, quasi conformally
flat or conformally symmetric \textit{LP}-Sasakian manifold endowed
with the structure$\ (\phi ,\xi ,\eta ,g)$ and $\gamma :I\rightarrow
M$ be a curve parametrized by arclength. Suppose that $\gamma $ is a
timelike curve that is its
velocity vector $T=\gamma ^{\prime }(s)$ is timelike. Let $%
\{T,N,B_{1},B_{2}\}$ be an orthonormal Frenet frame field tangent to $M$ along $%
\gamma $, where $N$ is the unit spacelike vector field in the direction $%
\nabla _{T}T$ , $B_{1}$ and $B_{2}$ are unit spacelike vectors. Then
by using the Frenet formulas (2.2.9), we have:
\begin{eqnarray*}
\tau _{2}(\gamma ) &=&\nabla _{T}^{3}T-R(T,\nabla _{T}T)T \\
&=&\nabla _{T}^{3}T-R(T,k_{1}N)T \\
&=&(3k_{1}k_{1}^{\prime })T+(k_{1}^{\prime \prime
}+k_{1}^{3}-k_{1}k_{2}^{2})N \\
&&+(2k_{1}^{\prime }k_{2}+k_{1}k_{2}^{\prime
})B_{1}+(k_{1}k_{2}k_{3})B_{2}-k_{1}R(T,N)T \\
&=&(3k_{1}k_{1}^{\prime })T+(k_{1}^{\prime \prime
}+k_{1}^{3}-k_{1}k_{2}^{2}-k_{1})N\\
&&+(2k_{1}^{\prime }k_{2}+k_{1}k_{2}^{\prime
})B_{1}+(k_{1}k_{2}k_{3})B_{2} \\
&=&0.
\end{eqnarray*}

It follows that $\gamma $ is a biharmonic curve if and only if%
\begin{eqnarray*}
k_{1}k_{1}^{\prime } &=&0 ,\\
k_{1}^{\prime \prime }+k_{1}(k_{1}^{2}-k_{2}^{2}-1) &=&0, \\
2k_{1}^{\prime }k_{2}+k_{1}k_{2}^{\prime } &=&0 ,\\
k_{1}k_{2}k_{3} &=&0.
\end{eqnarray*}%
If we look for nongeodesic solutions , that is for biharmonic curves with $%
k_{1}\neq 0,$ we obtain%
\begin{eqnarray*}
k_{1} &=&cons\tan t\neq 0,k_{2}=cons\tan t ,\\
k_{1}^{2}-k_{2}^{2} &=&1, \\
k_{2}k_{3} &=&0.
\end{eqnarray*}%
\end{proof}

\begin{proposition}
Let $\gamma :I\rightarrow S_{1}^{4}$ be a timelike nongeodesic
biharmonic curve parametrized by arclenght. Then
\begin{eqnarray}
\gamma^{(IV)}-2\gamma''+(1-k_{1}^{2})\gamma=0.
\end{eqnarray}
\end{proposition}

If $k_{1}=1$, then the general solution of (3.4) is
\begin{eqnarray*}
\gamma(s)=c_{1}+c_{2}s+c_{3}e^{-\sqrt{2}s}+c_{4}e^{\sqrt{2}s}
\end{eqnarray*}
Here $c_{i}$, $i=1,\,2,\,3,\,4$, are constant vectors. Since
$<\gamma(s),\gamma(s)>=1$ and $<\gamma'(s),\gamma'(s)>=-1$, by
choosing
\begin{eqnarray*}
c_{1}&=&(\frac{1}{\sqrt{2}},0,0,0,1),\,\,c_{2}=(0,0,0,0,0),\\
c_{3}&=&(-1,\frac{1}{\sqrt{2}},0,0,-\frac{1}{\sqrt{2}}),\,\,
c_{4}=(1,-\frac{\sqrt{2}}{4},\frac{1}{2\sqrt{2}},\frac{1}{2},\frac{1}{\sqrt{2}}),\,\,
\end{eqnarray*}
such that
\begin{eqnarray*}
<c_{1},c_{1}>&=&\frac{1}{2}\\
<c_{2},c_{2}>=<c_{3},c_{3}>&=&<c_{4},c_{4}>=0,\\
<c_{1},c_{2}>=<c_{1},c_{3}>&=&<c_{1},c_{4}>=0,\\
<c_{2},c_{3}>&=&<c_{2},c_{4}>=0,\\
<c_{3},c_{4}>&=&\frac{1}{4},
\end{eqnarray*}
we obtain following special solution of differential equation (3.4)
\begin{eqnarray*}
\gamma(s)&=&(\frac{1}{\sqrt{2}}-e^{-\sqrt{2}s}+e^{\sqrt{2}s},\frac{e^{-\sqrt{2}s}}{\sqrt{2}}-\frac{e^{\sqrt{2}s}}{2\sqrt{2}},\nonumber\\
&&\frac{e^{\sqrt{2}s}}{2\sqrt{2}},\frac{e^{\sqrt{2}s}}{2},1-\frac{e^{-\sqrt{2}s}}{\sqrt{2}}+\frac{e^{\sqrt{2}s}}{\sqrt{2}}),
\end{eqnarray*}
which is a circle.\\
If $k_{1}>1$, then the general solution of (3.4) is
\begin{eqnarray*}
\gamma(s)=c_{1}e^{as}+c_{2}e^{-as}+c_{3}\cos(bs)+c_{4}\sin(bs)
\end{eqnarray*} with $a=\sqrt{k_{1}+1}$ and $b=\sqrt{k_{1}-1}$. Here $c_{i}$, $i=1,\,2,\,3,\,4$, are constant
vectors. Since again $<\gamma(s),\gamma(s)>=1$ and
$<\gamma'(s),\gamma'(s)>=-1$, by choosing
\begin{eqnarray*}
c_{1}&=&(1,0,0,0,1),\,\,c_{2}=(-1,\frac{\sqrt{7}}{4},0,0,-\frac{3}{4}),\\
c_{3}&=&(0,0,\frac{1}{2},\frac{1}{2},0),\,\,
c_{4}=(-\frac{\sqrt{7}}{\sqrt{2}},\frac{1}{\sqrt{2}},0,0,-\frac{\sqrt{7}}{\sqrt{2}}),\,\,
\end{eqnarray*}
such that $c_{1}$ and $c_{2}$ are null vectors and
\begin{eqnarray*}
<c_{3},c_{3}>&=&<c_{4},c_{4}>=\frac{1}{b^{2}},\\
<c_{1},c_{2}>&=&\frac{1}{a^{2}},\\
<c_{1},c_{3}>&=&<c_{1},c_{4}>=0,\\
<c_{2},c_{3}>&=&<c_{2},c_{4}>=0,\\
<c_{3},c_{4}>&=&0,
\end{eqnarray*} with $a=2,\,\,b=\sqrt{2}$, we obtain following special
solution of differential equation (3.4)
\begin{eqnarray*}
\gamma(s)&=&(e^{2s}-e^{-2s}-\frac{\sqrt{7}}{\sqrt{2}}\sin(\sqrt{2}s),\frac{\sqrt{7}}{4}e^{-2s}+\frac{1}{\sqrt{2}}\sin(\sqrt{2}s),\nonumber\\
&&\frac{1}{2}\cos(\sqrt{2}s),\frac{1}{2}\cos(\sqrt{2}s),e^{2s}-\frac{3}{4}e^{-2s}-\frac{\sqrt{7}}{\sqrt{2}}\sin(\sqrt{2}s)),
\end{eqnarray*}
which is a helix with $k_{1}=3$ and $k_{2}=2\sqrt{2}$.

\begin{rem}
In this paper we do not consider the null curves in a
\textit{LP}-Sasakian manifold. Because a null curve in a
semi-Riemannian manifold can be considered as a 1-dimensional
degenerate submanifold and some difficulties arise when the
Laplacian operator is being defined in a degenerate submanifold.
Hence the biharmonicity of a null curve thought of a 1-dimensional
submanifold can not be defined by means of the variational problem.
\end{rem}

\textit{Authors' adress}:\\
Sad\i k KELE\c{S}, Selcen Y\"{U}KSEL PERKTA\c{S} and Erol
KILI\c{C},\\
Department of Mathematics,\\
Inonu University, 44280, Malatya/TURKEY\\
E-mail: skeles@inonu.edu.tr, selcenyuksel@inonu.edu.tr,
 ekilic@inonu.edu.tr
\end{document}